\documentclass[reqno,10pt]{amsart}
\usepackage[english]{babel}
\usepackage{amsthm}
\usepackage{amsmath}
\usepackage{amssymb}
\usepackage{amsfonts}
\usepackage{mathrsfs}
\usepackage{todonotes}
\usepackage{comment}

\usepackage[utf8]{inputenc}
\usepackage{latexsym}
\usepackage{verbatim}

\usepackage[all]{xy}
\usepackage{dsfont}
\usepackage[bookmarksnumbered,colorlinks,urlcolor=blue]{hyperref}
\usepackage{graphicx}
\usepackage{graphics}
\usepackage{float}
\usepackage{enumerate}
\def\bb#1\eb{\textcolor{blue}
{#1}} %
\def\br#1\er{\textcolor{red}
{#1}} %
\newcommand{\R}{\mathds R}

   \def\br#1\er{\textcolor{red}{#1}} %
      \def\bb#1\eb{\textcolor{blue}{#1}} %

\title[]{Locally extremal geodesic loops on Riemannian manifold}

\author[J.L. Flores]{Jos\'e Luis Flores}
\address{Departamento de \'Algebra, Geometr\'{\i}a y Topolog\'{\i}a,  Universidad de M\'alaga
\hfill\break\indent
Facultad de Ciencias, Campus Universitario de Teatinos,
\hfill\break\indent 29071 M\'alaga, Spain}
\email{floresj@uma.es}

\thanks{2010 {\em Mathematics Subject Classification:} Primary  53C22 \\
\textbf{Key words:} Geodesic loop, closed geodesic, conjugate point.}

\begin{document}
\newtheorem{thm}{Theorem}[section]
\newtheorem{prop}[thm]{Proposition}
\newtheorem{lemma}[thm]{Lemma}
\newtheorem{cor}[thm]{Corollary}
\theoremstyle{definition}
\newtheorem{defi}[thm]{Definition}
\newtheorem{notation}[thm]{Notation}
\newtheorem{exe}[thm]{Example}
\newtheorem{conj}[thm]{Conjecture}
\newtheorem{prob}[thm]{Problem}
\newtheorem{rem}[thm]{Remark}
\newtheorem{conv}[thm]{Convention}
\newtheorem{crit}[thm]{Criterion}
\newtheorem{claim}[thm]{Claim}

\newcommand{\ben}{\begin{enumerate}}
\newcommand{\een}{\end{enumerate}}

\newcommand{\bit}{\begin{itemize}}
\newcommand{\eit}{\end{itemize}}

\begin{abstract}
This note proves that any locally extremal non-self-conjugate geodesic loop in a Riemannian manifold is a closed geodesic. As a consequence, any complete and non-contractible Riemannian manifold with diverging injectivity radii along diverging sequences and without points conjugate to themselves, possesses a minimizing closed geodesic.
\end{abstract}

\maketitle

\vspace*{-.5cm}

\section{Introduction}\label{section1}

It is obvious that the classical result on the existence of a closed geodesic in any compact Riemannian manifold
is no longer true if we remove the compactness hypothesis. Consequently, any existence result of closed geodesic in the non-compact case require to assume further conditions on the topology of the manifold and on its geometry at infinity. In the last decades several significative theorems have appeared in this direction.
In \cite{Tho} Thorbergsson proved that any complete Riemannian manifold contains a closed geodesic provided that it is non-contractible and its sectional curvature is non-negative outside a compact set. In \cite{BG91, BG92}, Benci and Giannoni proved the existence of a closed geodesic in any complete Riemannian manifold with non-positive sectional curvature at infinity and whose free loop space has non-trivial homology in sufficiently high dimension. The proof involves variational techniques, and include both, a penalization method to overcome the possible noncompactness of $M$ and an estimate of the index form. More recently, in \cite{BGS} the authors extended previous result to the case of a complete Riemannian manifold {\em with boundary}, provided that the boundary is smooth, compact and {\em convex} (in the sense that its second fundamental form is non-negative in the direction of the inner normal).

The purpose of this note is to prove a characterization of closed geodesics in terms of locally extremal geodesics loops (Theorem \ref{t}), and deduce from it a new result about the existence of a minimizing closed geodesic in a (non-necessarily compact) complete Riemannian manifold with further conditions on the topology (non-contractible) and the geometry at infinity (diverging injectivity radii along diverging sequences), Theorem \ref{tt}.
%

In order to formalize these ideas, let us begin with some basic preliminaries.

Let $(M,g)$ be a complete Riemannian manifold. A (non-constant smooth) curve $\gamma:[0,l]\rightarrow M$ is a {\em geodesic loop} if $\gamma\mid_{(0,l)}$ is a unitary geodesic in $(M,g)$ and $\gamma(0)=\gamma(l)$. If, in addition,
$\dot{\gamma}(0)=\dot{\gamma}(l)$, it is a {\em closed (or periodic) geodesic}. A geodesic loop $\gamma:[0,l]\rightarrow M$ is said {\em locally minimizing} (resp. {\em locally maximizing}) if ${\rm length}(\gamma)\leq {\rm length}(\gamma')$ (resp. ${\rm length}(\gamma)\geq {\rm length}(\gamma')$) for any geodesic loop $\gamma':[0,l']\rightarrow M$ with initial conditions
($\gamma'(0)$, $\dot{\gamma}'(0))\in TM$ close\footnote{Here, $TM$ is implicitly endowed with, say, the Sasaki metric associated to $g$.} enough to the ones of $\gamma$, ($\gamma(0)$, $\dot{\gamma}(0)$). In general, a geodesic loop is said {\em locally extremal} if it is either locally minimizing or locally maximizing. Finally, a geodesic loop $\gamma:[0,l]\rightarrow M$ is said {\em self-conjugate} if $\gamma(0)$, $\gamma(l)$ are conjugate between them along $\gamma$.

\smallskip

This is the first important result of this note:

\begin{thm}\label{t}
Any locally extremal non-self-conjugate geodesic loop in a Riemannian manifold is a closed geodesic.
\end{thm}
%
%
%

\noindent The proof of Theorem \ref{t} will be developed in Section \ref{s2}, and is based on a length-shortening/lenghthening argument directly applied on geodesic loops by taking advantage of the exponential map and the absence of conjugate points. In the rest of this section we are going to apply this theorem to deduce a result about the existence of a minimizing closed geodesic in this ambient (Theorem \ref{tt}). To this aim, we begin by recalling the following well-known property (see \cite[Th. 13.3, p.239]{M}, \cite{Se}):

\begin{prop} For every point $p$ in a complete and non-contractible Riemannian manifold $(M,g)$ there exists some geodesic loop $\gamma:[0,l]\rightarrow M$ satisfying $\gamma(0)=\gamma(l)=p$.
\end{prop}
\noindent Next, assume that $(M,g)$ is a complete and non-contractible Riemannian manifold with {\em diverging injectivity radii along diverging sequences}. Then, the infimum of the injectivity radii of the points along $M$ is positive. In particular, the lengths of all geodesic loops on $M$ are bounded below by that positive number. Let $l>0$ be the infimum of such lengths and $\{\gamma_n\}\subset M$ be a sequence of geodesic loops in $(M,g)$ realizing that infimum, that is, ${\rm length}(\gamma_n)\searrow l$.
From the hypothesis about the divergence of the injectivity radii, the sequence $\{\gamma_n\}$ must remain in a bounded region of $M$, and so, up to a subsequence, $\dot{\gamma}_n(0)\rightarrow v$ for some unitary $v\in T_p M$, $p\in M$. Then, by continuity, one deduces that the unitary geodesic $\gamma$ with initial conditions $\gamma(0)=p$ and $\dot{\gamma}(0)=v$ satisfies $\gamma(l)=p$, and so, it is a minimizing (thus, simple) geodesic loop in $(M,g)$. In conclusion:
\begin{prop}\label{pp} Any complete and non-contractible Riemannian manifold with diverging injectivity radii along diverging sequences possesses a minimizing geodesic loop.
\end{prop}
\noindent This proposition joined to Theorem \ref{t} immediately provides the following result about existence of a minimizing closed geodesic on non-necessarily compact Riemannian manifolds under simple geometric and topological hypotheses:
\begin{thm}\label{tt} Any complete and non-contractible Riemannian manifold $(M,g)$ with diverging injectivity radii along diverging sequences possesses
\begin{itemize}
\item[(i)] either a minimizing self-conjugate geodesic loop,
\item[(ii)] or a minimizing closed geodesic.
\end{itemize}
In particular, if $M$ has no points conjugate to themselves, possibility (i) is ruled out, and the existence of a minimizing closed geodesic is ensured.
\end{thm}

\noindent  Note that Theorem \ref{tt} includes the classical result (cited at the beginning of this note) that any compact (thus, complete and non-contractible) Riemannian manifold admits some minimizing closed geodesic, under the additional hypothesis of non-existence of conjugate points to themselves.

\begin{rem} (1) In principle, the approach followed in this paper does not permit to circumvent the hypothesis about non-self-conjugacy of the locally extremal geodesic loop in Theorem \ref{t}. However, we suspect that this hypothesis is not really necessary. Consequently, we believe that assertion (i), and so, the hypothesis about non-existence of points conjugate to themselves can be removed from both, Theorem \ref{tt} and the paragraph above, resp.

(2) The pseudosphere shows that Proposition \ref{pp} and Theorem \ref{tt} are no longer true if we remove the hypothesis about the divergence of the injectivity radii along diverging sequences.
\end{rem}


\section{Proof of  Theorem \ref{t}} \label{s2}

Assume by contradiction that $\gamma:[0,l]\rightarrow M$ is a, say, locally minimizing, non-self-conjugate geodesic loop in a Riemannian manifold $(M,g)$ which is not a closed geodesic; that is, $\gamma:(0,l)\rightarrow M$ is a unitary geodesic, $\gamma(0)=\gamma(l)=p$, $\dot{\gamma}(0)\neq\dot{\gamma}(l)$ and $\gamma(0)(=\gamma(l))$ is not conjugate to itself along $\gamma$.
Then, $(d\exp_p)_{v_0}:T_{v_0}(T_p M)\rightarrow T_{\exp_{p}(v_0)} M$ is non-singular, where $v_0:=l\cdot\dot{\gamma}(0)$, and thus, $\exp_p(v_0)=\exp_{\gamma(0)}(l\cdot\dot{\gamma}(0))=\gamma(l)=p$.
By continuity, we can find $\delta>0$ small enough such that
$(d\exp_q)_{v}:T_{v}(T_q M)\rightarrow T_{\exp_q(v)} M$ is also non-singular, with $q:=\gamma(\delta)\neq p$, $v:=(l-\delta)\cdot\dot{\gamma}(\delta)$, and thus
\begin{equation}\label{d}
\exp_q(v)=\exp_{\gamma(\delta)}((l-\delta)\cdot\dot{\gamma}(\delta))=\gamma(\delta+(l-\delta))=\gamma(l)=p.
\end{equation}
So, by the Inverse Function Theorem, $\exp_q:U\subset T_q M\rightarrow V\subset M$
is a diffeomorphism between certain neighborhoods $U$ of $v\in T_q M$ and $V$ of $p\in M$. Even more, by taking $\delta>0$ smaller if necessary, we can assume that $V$ is a normal ball of center $p$ and radius $r_0$ greater than $\delta$. In particular, the geodesic segment $\alpha\equiv\gamma\mid_{[0,\delta]}$ is totally contained in $V$.
%
%
%
%
%
%

Consider the curve $$\overline{\alpha}:=(\exp_q\mid_{U})^{-1}\circ\alpha:[0,\delta]\rightarrow U\subset T_q M.$$
From (\ref{d}),
\[
v=(\exp_q\mid_{U})^{-1}(p)=(\exp_q\mid_{U})^{-1}(\alpha(0))=\overline{\alpha}(0)\in U\subset T_q M.
\]
So, if we define
\[
u:=\overline{\alpha}(\delta)=(\exp_q\mid_{U})^{-1}(\alpha(\delta))\in U\subset T_q M,
\]
the vectors $u$, $v$ are not only different (recall that $\alpha(0)=p\neq q=\alpha(\delta)$ and $\exp_q\mid_{U}$ is a diffeomorphism), but they also satisfy
\begin{equation}\label{m}
u/|u|\neq v/|v|.
\end{equation}
In fact, otherwise, $u=\lambda v$, for some positive $\lambda\neq 1$. If $0<\lambda<1$ then $\gamma\mid_{[\delta,\delta+\lambda(l-\delta)]}$ would be a geodesic loop,
\[
\gamma(\delta)=q=\exp_q(u)=\exp_q(\lambda v)=\exp_{\gamma(\delta)}(\lambda(l-\delta)\dot{\gamma}(\delta))=\gamma(\delta+\lambda(l-\delta)),
\]
with $(\gamma(\delta)=q,\dot{\gamma}(\delta))$ close to $(\gamma(0)=p,\dot{\gamma}(0))$, such that
\[
{\rm length}(\gamma\mid_{[\delta,\delta+\lambda(l-\delta)]})=\lambda(l-\delta)<l,
\]
in contradiction with the locally minimizing character of the geodesic loop $\gamma$. So, we can assume that $\lambda>1$. But, in this case, $\gamma\mid_{[l,\delta+\lambda(l-\delta)]}$ is a geodesic segment connecting $\gamma(l)=p$ with $\gamma(\delta+\lambda(l-\delta))=\exp_{\gamma(\delta)}(\lambda(l-\delta)\dot{\gamma}(\delta))=\exp_q(\lambda v)=\exp_q(u)=q$, and satisfying
\[
{\rm length}(\gamma\mid_{[l,\delta+\lambda(l-\delta)]})<\lambda(l-\delta)=\lambda |v|=|u|<r_0.
\]
So, $\gamma\mid_{[l,\delta+\lambda(l-\delta)]}$ is totally contained in the normal ball $V$, and thus, it must coincide with $\gamma\mid_{[0,\delta]}$, in contradiction with $\dot{\gamma}(0)\neq\dot{\gamma}(l)$.

Next, observe that
\[
\alpha(t)=\exp_q(\overline{\alpha}(t))=f(r(t),t),\quad t\in [0,\delta],
\]
where
\begin{equation}\label{mm}
f(r,t):=\exp_q(r\cdot w(t)),\quad r(t):=|\overline{\alpha}(t)|,\quad w(t):=\frac{\overline{\alpha}(t)}{|\overline{\alpha}(t)|}.
\end{equation}
In particular,
\[
\dot{\alpha}(t)=\frac{d}{dt}f(r(t),t)=\frac{\partial f}{\partial r}\dot{r}(t)+\frac{\partial f}{\partial t}.
\]
From the Gauss Lemma, $g(\partial f/\partial r,\partial f/\partial t)=0$. Moreover, $|\partial f/\partial r|=1$.
Hence,
\[
|\dot{\alpha}(t)|^2=\left|\frac{\partial f}{\partial r}\right|^2 \dot{r}(t)^2+\left|\frac{\partial f}{\partial t}\right|^2
=\dot{r}(t)^2+\left|\frac{\partial f}{\partial t}\right|^2,
\]
and thus,
\[
|\dot{\alpha}(t)|=\sqrt{\dot{r}(t)^2+\left|\frac{\partial f}{\partial t}\right|^2}=|\dot{r}(t)|+\frac{\left|\frac{\partial f}{\partial t}\right|^2}{\sqrt{\dot{r}(t)^2+\left|\frac{\partial f}{\partial t}\right|^2}+|\dot{r}(t)|}.
\]
Therefore,
\begin{equation}\label{eqq}
\begin{array}{l}
{\rm length}(\alpha)=\int_0^{\delta}|\dot{\alpha}(t)| dt\geq\int_0^{\delta}\dot{r}(t)dt + \int_{0}^{\delta}\frac{\left|\frac{\partial f}{\partial t}\right|^2}{\sqrt{\dot{r}(t)^2+\left|\frac{\partial f}{\partial t}\right|^2}+|\dot{r}(t)|}dt
\\
\qquad\qquad\qquad\qquad =|\overline{\alpha}(\delta)|-|\overline{\alpha}(0)|+\int_{0}^{\delta}\frac{\left|\frac{\partial f}{\partial t}\right|^2}{\sqrt{\dot{r}(t)^2+\left|\frac{\partial f}{\partial t}\right|^2}+|\dot{r}(t)|}dt.
\end{array}
\end{equation}
From (\ref{m}),
\[
w(\delta)=\frac{\overline{\alpha}(\delta)}{|\overline{\alpha}(\delta)|}=\frac{u}{|u|}\neq \frac{v}{|v|}=\frac{\overline{\alpha}(0)}{|\overline{\alpha}(0)|}=w(0).
\]
Hence, $|\partial f/\partial t|>0$ for some $t\in [0,\delta]$ (recall (\ref{mm}) and the fact that $\exp_q\mid_U$ is a diffeomorphism). In particular,
\begin{equation}\label{as}
\int_{0}^{\delta}\frac{\left|\frac{\partial f}{\partial t}\right|^2}{\sqrt{\dot{r}(t)^2+\left|\frac{\partial f}{\partial t}\right|^2}+|\dot{r}(t)|}dt=\eta_0>0.
\end{equation}
Moreover, we know that
\begin{equation}\label{a}
|\overline{\alpha}(0)|=|v|=|(l-\delta)\dot{\gamma}(\delta)|=l-\delta,\qquad |\overline{\alpha}(\delta)|=|u|.
\end{equation}
So, from (\ref{eqq}), (\ref{as}), (\ref{a}),
\[
\delta={\rm length}(\alpha)\geq |u|-(l-\delta)+\eta_0,\quad\hbox{and thus, $\quad |u|<l={\rm length}(\gamma)$.}
\]
Summarizing, $\gamma':[0,|u|]\rightarrow M$, $\gamma'(s)=\exp_q(s\cdot u/|u|)$, is a geodesic loop,
\[
\gamma'(0)=\exp_q(0)=q=\alpha(\delta)=\exp_q(u)=\gamma'(|u|),
\]
with $(\gamma'(0)=q,\dot{\gamma}'(0)=u/|u|)$ close to $(\gamma(0)=p,\dot{\gamma}(0))$, such that
\[
{\rm length}(\gamma')=|u|<{\rm length}(\gamma),
\]
in contradiction with the locally minimizing character of $\gamma$.

Finally, if, instead of locally minimizing, we initially assume that the geodesic loop $\gamma:[0,l]\rightarrow M$ is locally maximizing, we can arrive to a contradiction again by taking $\delta$ negative instead of positive in previous argument.

\section*{Acknowledgments}
The author is partially supported by the project MTM2016-78807-C2-2-P (Spanish MINECO with FEDER funds). He thanks Miguel S\'anchez for valuable suggestions.

\end{document}